\documentclass[11pt]{article}
\newtheorem{theorem}{Theorem}[section]
\newtheorem{lemma}[theorem]{Lemma}

\newtheorem{conjecture}[theorem]{Conjecture}

\usepackage{color}
\usepackage{amssymb,amsmath}
\usepackage{mathrsfs}
\usepackage{hyperref}
\DeclareMathOperator{\ri}{ri}
\DeclareMathOperator{\sgn}{sgn}
\def\f{\noindent}
\def\h{\hfill  $\Box$\vspace{10pt}}
\def\p{\f {\bf Proof}\hskip10pt}

\begin{document}

\title{Sign patterns that require $\mathbb{H}_n$ exist for each $n\geq 4$ }
\author{Wei Gao\thanks{Corresponding author. E-mail: wgao2@gsu.edu \vspace{6pt}}, Zhongshan Li, Lihua Zhang\\
Dept of Math and Stat, Georgia State University, Atlanta, GA 30302-4110, USA}
\date{}
\maketitle
\begin{abstract}
The refined inertia of a square real matrix $A$ is the ordered $4$-tuple $(n_+, n_-, n_z, 2n_p)$,
where $n_+$ (resp., $n_-$) is the number of eigenvalues of $A$ with positive (resp., negative) real part,
$n_z$ is the number of zero eigenvalues of $A$,
and $2n_p$ is the number of nonzero pure imaginary eigenvalues of $A$.
For $n \geq 3$, the set of refined inertias $\mathbb{H}_n=\{(0, n, 0, 0), (0, n-2, 0, 2), (2, n-2, 0, 0)\}$ is important for the onset of Hopf bifurcation in
dynamical systems.
We say that an $n\times n$ sign pattern ${\cal A}$  requires $\mathbb{H}_n$ if
$\mathbb{H}_n=\{\ri(B) | B \in Q({\cal A})\}$.
Bodine et al.  conjectured that no $n\times n$ irreducible sign pattern that requires $\mathbb{H}_n$ exists for $n$ sufficiently large, possibly $n\ge 8$.
However, for each $n \geq 4$, we identify three $n\times n$ irreducible sign patterns that require $\mathbb{H}_n$, which  resolves this conjecture.

\bigskip
\f {\it AMS classification:}  15B35, 15A18, 05C50

\f {\bf Keywords:} Eigenvalues; Refined inertia; Sign pattern.
\end{abstract}

%%%%%%%%%%%%%%% Section 1
\section{Introduction}
\hskip\parindent
The refined inertia of a square real matrix $B$, denoted $\ri(B)$, is the ordered $4$-tuple
$(n_+(B),$ $ n_-(B), n_z(B), 2n_p(B))$,
where $n_+(B)$ (resp., $n_-(B)$) is the number of eigenvalues of $B$ with positive
(resp., negative) real part,
$n_z(B)$ is the number of zero eigenvalues of $B$,
and $2n_p(B)$ is the number of pure imaginary eigenvalues of $B$.  Refined inertias were introduced in \cite{Deaett10}, and have been the focus of study in several recent papers such as \cite{Bodine12,Gao14,Garnett13,Garnett14,Olesky13,Yu12}.

A {\it sign pattern} ({\it matrix}) is a matrix whose entries are from the set $\{+, -, 0\}$. For a real matrix $B$, sgn($B$) is the sign pattern matrix obtained by replacing each positive (resp., negative, zero) entry of $B$ by $+$ (resp., $-, 0$). For an $n\times n$ sign pattern matrix $\cal A$, the qualitative class
of $\cal A$, denoted $Q(\cal A)$, is defined as
$Q({\cal A})=\{B\in M_n({\mathbb{R}}) \mid \hbox{sgn}(B)={\cal A}\}$.

A {\it permutation sign pattern} is a square sign pattern with entries from the set $\{0, +\}$, where the entry $+$ occurs precisely once in each row and in each column. A {\it signature sign pattern} is a square diagonal sign pattern all of whose diagonal entries are nonzero. Let ${\cal A}_1$ and ${\cal A}_2$ be two square sign patterns of the same order.
Sign pattern ${\cal A}_1$ is said to be {\it permutationally similar} to ${\cal A}_2$ if there exists a permutation sign pattern ${\cal P}$ such that ${\cal A}_2 ={\cal P}^{T} {\cal A}_1 {\cal P}$.
Sign pattern ${\cal A}_1$ is said to be {\it signature similar} to ${\cal A}_2$ if there exists a signature sign pattern ${\cal D}$ such that ${\cal A}_2 ={\cal D}{\cal A}_1 {\cal D}$.
Two sign patterns are said to be \emph{equivalent} if one can be obtained from the other by transposition,
signature similarity, permutation similarity, or any combination of these.

Let $n\ge 3$ and let $\mathbb{H}_n=\{(0, n, 0, 0),\ (0, n-2, 0, 2),\ (2, n-2, 0, 0)\}$. As pointed out by Bodine et al.~\cite{Bodine12}, $\mathbb{H}_n$
is an important set of refined inertias which can signal the onset of periodic solutions by Hopf bifurcation in dynamical systems.
We say that an $n\times n$ sign pattern ${\cal A}$  {\em requires} $\mathbb{H}_n$ if
$\mathbb{H}_n=\{\ri(B) | B \in Q({\cal A})\}$,
and  ${\cal A}$ {\em allows} $\mathbb{H}_n$ if
$\mathbb{H}_n\subseteq \{\ri(B) | B \in Q({\cal A})\}$.
In \cite{Bodine12}, the authors made the following conjecture.

\begin{conjecture}[\cite{Bodine12}]\label{conj}\
No $n\times n$ irreducible sign pattern that requires $\mathbb{H}_n$ exists for $n$ sufficiently large, possibly $n\ge 8$.
\end{conjecture}

In this paper, for each $n\geq 4$, we identify three $n\times n$ irreducible sign patterns that each require $\mathbb{H}_n$,
which negatively resolves the preceding conjecture.

Let ${\cal A}_1$, ${\cal A}_2$, ${\cal A}_3$ be sign patterns of order $n \geq 4$ defined by
$${\cal A}_1=\left[\begin{array}{cccccc}
+&+&+&\cdots&+&+\\
-&0&\\
-&&-&\\
\vdots&&&\ddots\\
-&&&&-&\\
-&&&&&-
\end{array}\right],\ \ \
{\cal A}_2= \left[\begin{array}{cccccc}
-&+&+&\cdots&+&+\\
-&0&\\
+&&-&\\
\vdots&&&\ddots\\
+&&&&-&\\
+&&&&&-
\end{array}\right],
$$
$$
{\cal A}_3=\left[\begin{array}{cccccc}
-&+&+&\cdots&+&+\\
+&0&\\
+&&-&\\
\vdots&&&\ddots\\
+&&&&-&\\
-&&&&&+
\end{array}\right],
$$
where  all the off-diagonal entries except those on the first row or first column are zeros.
We will show that for $i=1,2,3$, the ${\cal A}_i$ require $\mathbb{H}_n$ for each $n\ge 4$.

Throughout what follows, we let $B$ denote a real matrix of order $n \geq 4$ of the form
$$B=\left[\begin{array}{cccccccc}
a_1&1&1&\cdots&1&1\\
a_2&0&&\\
a_3&&-b_1&\\
\vdots&&&\ddots\\
a_{n-1}&&&&-b_{n-3}\\
a_{n}&&&&&-b_{n-2}
\end{array}\right],\eqno{(1.1)}$$
where $b_j>0$ for $j=1,2,\dots, n-3$, and
suitable real values for $a_1, a_2, \dots, a_n$ and $b_{n-2}$ are taken so that
$B\in Q({\cal A}_i)$ for some $i \in \{ 1,2,3 \}$.

%%%%%%%%%%%%% Section 2
\section{The sign patterns ${\cal A}_i$ allow  $\mathbb{H}_n$ for $i=1,2,3$}
\hskip\parindent
Note that the $4\times 4$ sign patterns ${\cal A}_1$, ${\cal A}_2$ and ${\cal A}_3$ given in Section 1 are equivalent to ${\cal S}_5$,
${\cal S}_2$ and ${\cal S}_4$ defined in \cite[p.624]{Garnett13}, respectively.
Thus the next lemma follows from Theorems 2.5, 2.6, and 2.8 in \cite{Garnett13}.

\begin{lemma}[\cite{Garnett13}]\label{lemma:H4}\
The $4\times 4$ sign patterns ${\cal A}_i$ require  $\mathbb{H}_4$ for $i=1,2,3$.
\end{lemma}

We now show that for $i=1,2,3$, the $n\times n$ sign patterns ${\cal A}_i$ allow $\mathbb{H}_n$ for each $n\ge 5$.

\begin{theorem}\label{allow:Hn}\
For $i=1,2,3$, the sign patterns ${\cal A}_i$ allow $\mathbb{H}_n$ for each $n\ge 5$, .
\end{theorem}

\p
Choose any $i \in \{1, 2, 3\}$. By Lemma \ref{lemma:H4}, for each of the refined inertias
$(0,4,0,0)$, $(0,2,0,2)$ and $(2,2,0,0)$,
there exist  suitable values of $a_1,a_2,a_3,a_4,b_1,b_2$ such that
$$B_{4\times 4}=\left[\begin{array}{cccc}
a_1&1&1&1\\
a_2&0&\\
a_3&&-b_1&\\
a_4&&&-b_2
\end{array}\right]\in Q({\cal A}_i)$$
has this refined inertia.

For $n\ge 5$, consider the $n\times n$ matrix
$$B=\left[\begin{array}{ccccccc}
a_1&1&1&1&\cdots&1&1\\
a_2&0&\\
\frac{a_3}{n-3}&&-b_1&\\
\frac{a_3}{n-3}&&&-b_1&\\
\vdots&&&&\ddots\\
\frac{a_3}{n-3}&&&&&-b_1&\\
a_4&&&&&&-b_2
\end{array}\right].
$$
Then $B\in Q({\cal A}_i)$, and
$$|\lambda I-B|=\left|\begin{array}{ccccccc}
\lambda-a_1 &-1&-1&-1&\cdots&-1&-1\\
-a_2&\lambda &\\
-\frac{a_3}{n-3}&&\lambda +b_1&\\
-\frac{a_3}{n-3}&&&\lambda +b_1&\\
\vdots&&&&\ddots\\
-\frac{a_3}{n-3}&&&&&\lambda +b_1&\\
-a_4&&&&&&\lambda +b_2
\end{array}\right|$$
$$=\left|\begin{array}{ccccccc}
\lambda-a_1 &-1&-1&-1&\cdots&-1&-1\\
-a_2&\lambda &\\
-a_3&&\lambda +b_1&\lambda +b_1&\cdots&\lambda +b_1\\
-\frac{a_3}{n-3}&&&\lambda +b_1&\\
\vdots&&&&\ddots\\
-\frac{a_3}{n-3}&&&&&\lambda +b_1&\\
-a_4&&&&&&\lambda +b_2
\end{array}\right|$$
$$=\left|\begin{array}{ccccccc}
\lambda-a_1 &-1&-1&0&\cdots&0&-1\\
-a_2&\lambda &\\
-a_3&&\lambda +b_1\\
-\frac{a_3}{n-3}&&&\lambda +b_1&\\
\vdots&&&&\ddots\\
-\frac{a_3}{n-3}&&&&&\lambda +b_1&\\
-a_4&&&&&&\lambda +b_2
\end{array}\right|$$
$$=(\lambda +b_1)^{n-4}\left|\begin{array}{ccccccc}
\lambda-a_1 &-1&-1&-1\\
-a_2&\lambda &\\
-a_3&&\lambda +b_1\\
-a_4&&&\lambda +b_2
\end{array}\right|
=(\lambda +b_1)^{n-4}|\lambda I-B_{4\times 4}|.$$
So the multiset of the eigenvalues of $B$ is given by  $\sigma(B)=\{-b_1,\dots,-b_1\}\cup \sigma(B_{4\times 4})$, in which each set is interpreted as a multiset.
It follows  that $n_-(B)=n_-(B_{4\times 4})+(n-4)$, $n_+(B)=n_+(B_{4\times 4})$,
$n_z(B)=n_z(B_{4\times 4})$, and $2n_p(B)=2n_p(B_{4\times 4})$.
Thus the $n\times n$ sign pattern ${\cal A}_i$ allows
$\mathbb{H}_n=\{(0,n,0,0),\  (0,n-2,0,2),\  (2,n-2,0,0)\}$. \hfill \h

%%%%%%%%%%%%%%%%%%% Section 3

\section{The main result}
\hskip\parindent
In this section, we establish that for $i=1,2,3$, the $n\times n$ sign patterns ${\cal A}_i$ require
$\mathbb{H}_n$ for each $n\ge 5$. As in the introduction, throughout this section we let $B$ be a
real matrix in the form $(1.1)$ of order $n \geq 5$.

First, we consider the case that all the $b_j$ are distinct  for $j=1, 2, \dots,n-2$ 
and show the following result.

\begin{theorem} \label{thm:distint_b_i}\
Let $n\ge 5$ and let $B$ have the form $(1.1)$.
If all the $b_j$ are distinct for $j=1, 2, \dots,n-2$, then $\ri(B) \in \mathbb{H}_n$.
\end{theorem}

Since all the $b_j$ are distinct for $j=1, 2, \dots,n-2$, we may assume, applying a permutation similarity if necessary, that
$B$ is subjected to
 $b_1>b_2>\dots>b_{n-2}$ in Theorem \ref{thm:distint_b_i}. To prove Theorem \ref{thm:distint_b_i}, we need the following lemmas. We also assume that $b_1>b_2>\dots>b_{n-2}$ in Lemmas \ref{lemma3.2}--\ref{lemma-even-odd}.

\begin{lemma}\label{lemma3.2}\
For $j=1, 2, \dots, n-2$,
$$|b_jI+B|=-a_{j+2}b_j\prod_{\substack{m=1 \\ m\ne j}}^{n-2}(b_j-b_m).$$
\end{lemma}

\p Note that row $j+2$ as well as column $j+2$ of $b_jI+B$ has exactly one nonzero entry, namely the first entry,
which may be used to zero out all other entries in the first row or the  first column without affecting the determinant.
Hence,
\begin{align}
|b_jI+B| & =\left| \begin{array}{ccccccccc}
0&0&0&\cdots&0&1&0&\cdots&0\\
0&b_j&\\
0&&b_j-b_1&\\
\vdots&&&\ddots\\
0&&&&b_j-b_{j-1}&\\
a_{j+2}&&&&&0\\
0&&&&&&b_j-b_{j+1}\\
\vdots&&&&&&&\ddots\\
0&&&&&&&&b_j-b_{n-2}
\end{array}\right|    \notag   \\
& =-a_{j+2}b_j\prod_{\substack{m=1\\ m\ne j}}^{n-2}(b_j-b_m).   \hspace{7.5cm} \Box \notag
\end{align}

In view of  Lemma \ref{lemma3.2}, the following two results are straightforward.

\begin{lemma}\label{lemma3.4}\
Suppose $B \in Q({\cal A}_i)$ for $i \in \{1,2,3\}$. Then
\begin{equation*}
\hbox{\em sgn}(\det(b_jI+B)) = \begin{cases}
(-)^{j+1} &\text{for $j=1,2,\dots,n-2$ if $i =1$};\\
(-)^{j} &\text{for $j=1,2,\dots,n-2$ if $i =2$};\\
(-)^{j} &\text{for $j=1,2,\dots,n-3$ if $i =3$}.
\end{cases}
\end{equation*}
\end{lemma}

\begin{lemma}\label{lemma3.3}\
The eigenvalues of $B$ do not include $-b_j$ for any $j \in \{1,2,\dots,n-2\}$.
\end{lemma}

\begin{lemma}  \label{remark3.5}\
Suppose $B \in Q({\cal A}_i)$ for $i \in \{1,2,3\}$. Then
$n_{-}(B)\geq n-4 \geq 1$. Furthermore,  if $i \neq 3$, $n_{-}(B) \geq n-3$.
\end{lemma}

\p
Observe that by Lemma \ref{lemma3.4}, the real function $p(t) = \det (tI -B)$ takes on nonzero values
of opposite signs at $-b_j$ and $-b_{j+1}$, for  $j=1, 2, \dots, n-4$.
Thus, by the Intermediate Value Theorem, $p(t)$ has at least one real zero in each open interval
$(-b_j,-b_{j+1})$. It follows that the matrix $B$ has at least one real eigenvalue in
$(-b_j,-b_{j+1})$, for $j=1, 2, \dots, n-4$. Thus $n_{-}(B)\geq n-4 \geq 1$.
Furthermore, if $i \in\{1, 2\}$, then by Lemma \ref{lemma3.4},
$B$ has at least one real eigenvalue in $(-b_j,-b_{j+1})$, for $j=1, 2, \dots, n-3$,
so we have  $n_{-}(B) \geq n-3$.
\h

\begin{lemma}  \label{remark3.6}
 $\sgn(\det(B))=(-)^n$. Furthermore, $n_z(B)=0$, and $n_-(B)$ and $n$ have the same parity. 
\end{lemma}

\p Expanding the determinant along the second column reveals that
$\hbox{sgn}(\det(B))=(-)^n$.  Consequently,  $n_z(B)=0$. It follows that  $\sgn(\det(B))=(-)^n = (-)^{n_-(B)}$. Hence,  
 $n_-(B)$ and $n$ have the same parity.
 \hfill \h

For any $r \in {\mathbb{R}}$, define $\Delta(r)$ to be the number of eigenvalues $\lambda$ of $B$
in the closed left half-plane with $\hbox{Re}(\lambda) \le -r$.
It is clear that
$$n_-(B)\ge\Delta(b_{n-3})=\Delta(b_1)
+\sum_{j=1}^{n-4}[\Delta(b_{j+1})-\Delta(b_j)].   \eqno{(3.1)}$$

\begin{lemma}\label{lemma3.6}\
For $j=1,2,\dots, n-3$,
$n_-(b_jI+B)$ and $\Delta(b_j)$ have the same parity.
\end{lemma}

\p
Note that $\lambda$ is an eigenvalue of $B$ if and only if
$b_j+\lambda$ is an eigenvalue of $b_jI+B$,
that the non-real eigenvalues of $b_jI+B$ occur in conjugate pairs,
and that $-b_1, -b_2, \dots, -b_{n-3}$ are not eigenvalues of $B$ by Lemma \ref{lemma3.3}.
We see that for $j=1,2,\dots, n-3$,
\begin{itemize}
\item 
$n_-(b_jI+B)=$ the number of eigenvalues $\lambda$ of $B$ satisfying $\hbox{Re}(\lambda)< -b_j$;
\item 
$\Delta(b_j)=$ the number of eigenvalues $\lambda$ of $B$ satisfying $\hbox{Re}(\lambda)\le -b_j$;
\item 
the number of eigenvalues $\lambda$ of $B$ satisfying $\hbox{Re}(\lambda)= -b_j$ is even.
\end{itemize}
So $n_-(b_jI+B)$ and $\Delta(b_j)$ have the same parity.
\hfill  \h

\begin{lemma}\label{lemma-even-odd}\
Let $k=n_-(B)$. Then $k\ge n-2$.
\end{lemma}

\p
If $i=1$ or $i=2$, then by Lemma \ref{remark3.5}, we have $k=n_-(B)\geq n-3$. By Lemma \ref{remark3.6},
$k$ and $n$ have the same parity. It follows that $k \geq n-2$, as desired. Hence, assume $i = 3$.

We claim that for every $j \le n-3$, the parity of $j$ and $\Delta(b_{j})$ are the same.
Otherwise, if there exists an even index  $j$  $ \le n-3 $ such that $\Delta(b_{j})$ is odd,
by Lemmas \ref{lemma3.4} and \ref{lemma3.6},
we have that $\det(b_{j}I+B)>0$ and $n_-(b_{j}I+B)$ is odd,
which is a contradiction; if there exists an odd index $j$ $ \le n-3 $ such that $\Delta(b_{j})$ is even,
by Lemmas \ref{lemma3.4} and \ref{lemma3.6},
we have that $\det(b_{j}I+B)<0$ and  $n_-(b_{j}I+B)$ is even, which is a contradiction.

Thus
$\Delta(b_{1})$ is odd, and
$\Delta(b_{j+1})-\Delta(b_{j}) > 0 $ is odd for $1\le j\le n-4$.
So by (3.1),
\begin{align}
k& \ge\Delta(b_{n-3})=\Delta(b_1)
+\sum_{j=1}^{n-4}[\Delta(b_{j+1})-\Delta(b_j)]
\ge n-3. \notag
\end{align}
By Lemma \ref{remark3.6}, $k$ and $n$ have the same parity.  It follows that $k \geq n-2$. \hfill \h

We now complete the proof of Theorem \ref{thm:distint_b_i}.

\vspace{2mm}
\noindent {\bf Proof  of Theorem 3.1}
By Lemmas \ref{remark3.6} and \ref{lemma-even-odd}, we have $n_z(B)=0$ and $n_-(B)=n-2$ or $n_-(B)=n$. 
It follows  that $\ri(B)\in \mathbb{H}_n$. \hfill  \h

We are now ready to establish the main result.

\begin{theorem}\label{thm:main}\
For $i=1,2,3$, the $n\times n$ sign patterns ${\cal A}_i$ require $\mathbb{H}_n$ for each $n\ge 4$.
\end{theorem}

\p  Fix any $i \in \{ 1, 2,3 \}$.  We proceed by induction on  the order $n$ of ${\cal A}_i$.

By Lemma \ref{lemma:H4}, the result holds for $n=4$.

Suppose that the $(n-1)\times (n-1)$ sign pattern ${\cal A}_i$  requires $\mathbb{H}_{n-1}$ for some $n\ge 5$.
We prove that the $n\times n$ sign pattern ${\cal A}_i$  requires $\mathbb{H}_n$.
By Theorem \ref{allow:Hn}, ${\cal A}_i$  allows $\mathbb{H}_n$.
Thus we only need to prove that $\ri(B) \in \mathbb{H}_n$ for every  $B \in Q({\cal A}_i)$.

For any $B \in Q({\cal A}_i)$, by performing a diagonal similarity on $B$ if necessary, we may assume that $B$ has the form (1.1).
If all the $b_j$ are distinct for $j=1, 2, \dots,n-2$, then by Theorem \ref{thm:distint_b_i} $\ri(B) \in \mathbb{H}_n$.

Now suppose that two of the  $b_j$ are the same for $j=1, 2, \dots,n-2$. Note that in the case of $B\in Q({\cal A}_3)$,
$b_{n-2}$ is different from each  $b_j$  with $j\leq n-3$ as $b_j>0>b_{n-2}$.
By performing a permutational similarity if necessary,
without loss of generality, we may assume that $b_1=b_2$. Then
$$|\lambda I-B|=\left|\begin{array}{ccccccc}
\lambda-a_1&-1&-1&-1&-1&\cdots&-1\\
-a_2&\lambda&\\
-a_3&&\lambda+b_1&\\
-a_4&&&\lambda+b_1&\\
-a_5&&&&\lambda+b_3\\
\vdots&&&&&\ddots\\
-a_n&&&&&&\lambda+b_{n-2}
\end{array}\right|$$
$$=\left|\begin{array}{ccccccc}
\lambda-a_1&-1&-1&-1&-1&\cdots&-1\\
-a_2&\lambda&\\
-a_3&&\lambda+b_1&\\
-a_3-a_4&&\lambda+b_1&\lambda+b_1&\\
-a_5&&&&\lambda+b_3\\
\vdots&&&&&\ddots\\
-a_n&&&&&&\lambda+b_{n-2}
\end{array}\right|$$
$$=\left|\begin{array}{ccccccc}
\lambda-a_1&-1&0&-1&-1&\cdots&-1\\
-a_2&\lambda&\\
-a_3&&\lambda+b_1&\\
-a_3-a_4&&&\lambda+b_1&\\
-a_5&&&&\lambda+b_3\\
\vdots&&&&&\ddots\\
-a_n&&&&&&\lambda+b_{n-2}
\end{array}\right|$$
$$=(\lambda+b_1)\left|\begin{array}{cccccc}
\lambda-a_1&-1&-1&-1&\cdots&-1\\
-a_2&\lambda&\\
-a_3-a_4&&\lambda+b_1&\\
-a_5&&&\lambda+b_3\\
\vdots&&&&\ddots\\
-a_n&&&&&\lambda+b_{n-2}
\end{array}\right|.$$
Take the $(n-1)\times (n-1)$ matrix
$$B_1=\left[\begin{array}{cccccc}
a_1&1&1&1&\cdots&1\\
a_2&0&\\
(a_3+a_4)&&-b_1\\
a_5&&&-b_3\\
\vdots&&&&\ddots\\
a_n&&&&&-b_{n-2}
\end{array}\right].$$
Then
$$\sigma(B)=\{-b_1\}\cup \sigma(B_1),$$
\mbox{where the sets are interpreted as multisets.}
Note that  $B_1 \in Q({\cal A}_i)$  has  order $n-1$.
By the induction  hypothesis, ${\cal A}_i$ of order $n-1$
requires $\mathbb{H}_{n-1}=\{(0,n-1,0,0),\ (0,n-3,0,2),\ (2,n-3,0,0)\}$.
Thus $\ri(B)$ is one of $(0, n ,0,0)$, $(0,n-2,0,2)$ and $(2,n-2,0,0)$.
It follows that  $\ri(B)\in \mathbb{H}_n$.

This completes the proof.   \hfill \h

\section*{Acknowledgement}
\hskip\parindent
The authors express their sincere thanks to the anonymous referee
for valuable suggestions which greatly improved the exposition of the paper.

\end{document}